\documentclass[11pt]{article}
\usepackage[top=2cm,bottom=2cm,left=2cm,right=2cm]{geometry}
\usepackage{amsmath,amsthm,amssymb}
\usepackage{graphicx,enumerate,url,color}

\numberwithin{equation}{section}

\begin{document}

\title{Calculating the day of the week: null-days algorithm}
\author{Xiang-Sheng Wang\footnote{Email: xswang@semo.edu}\\
Department of Mathematics\\Southeast Missouri State University\\Cape Girardeau, Missouri 63701}

\maketitle

\begin{abstract}
In this paper, we propose a new algorithm of calculating the day of the week for any given century, year, month and day in Gregorian calendar. We provide two simple formulas to convert the century and the year into two integers. Then we introduce a list of null-days to transform the month and the day into another integer. Adding these three integers together and calculating the sum's residue modulo $7$ gives a number between $0$ and $6$, which corresponds to Sunday until Saturday. Slight modification is needed if we have a leap year and the given month is either January or February. Our null-days algorithm is simple enough to be done by mental calculation, and the list of null-days has memorable patterns.
\end{abstract}

\section{Introduction}
It is interesting to see someone calculating the day of the week in mind. In this paper, we will propose a simple algorithm to fulfill this task with very basic mathematics and very few to be memorized.

For any given Gregorian calendar date, we can isolate the century item ($c$), year item ($y$), month item ($m$) and day item ($d$), respectively. Here, $m$ denotes the number of months in the year, i.e., $m=1$ for January, $m=2$ for February, and so on. Early in 1887, Lewis Carroll \cite{Ca87} devised a method to convert the items $c,y,m,d$ into a number $w$ (week), which modulo $7$ gives the day of the week. That is, $w\equiv1$ denotes Monday, $w\equiv2$ stands for Tuesday, and so on. The most complicated parts in Carroll's method are the translations of year item and month item, which involve a division of $12$ and additions of large numbers. To avoid addition of large numbers, one may assign certain numbers for the months, which is referred to as table method. The disadvantage of table method is that the months table is so irregular that it is easy to forget. John H. Conway \cite{Co73} introduced the famous doomsday rule to resolve this short back. Conway's doomsday consists of several memorable dates such as Valentine's Day, Boxing Day and the dates where $m=d=4,6,8,10,12$.

In this paper, we will modify Conway's doomsday rule and make a list of so called null-days which have certain memorable patterns. In addition, we provide two simple formulas to convert the century item and year item into two numbers. Coupling these two steps, we will be able to do mental calculation of the day of the week for any given date in Gregorian calendar.

\section{The null-days algorithm}
First, we introduce a list of null-days:
$$1/1;~~~ 3/5,~ 5/7,~ 7/9,~ 9/3;~~~ 2/12,~ 12/10,~ 10/8,~ 8/6,~ 6/4,~ 4/2;~~~ 11/12.$$
It is easy to verify that except for leap years, the above null-days in the same year share the same day of the week. The null-days are memorable in the sense that $1/1$ is just the first day of the year; the dates $3/5, 5/7, 7/9, 9/3$ are simply rotating the odd numbers $3,5,7,9$; the dates $2/12, 12/10, 10/8, 8/6, 6/4, 4/2$ are rotating the even numbers $12, 10, 8, 6, 4, 2$. The only exceptional null-day is $11/12$, which could be memorized in the way that $11$ (November) is equivalent with $1+1=2$ (February). We denote by $w_0(m,d)$ the difference of the date $m/d$ from the corresponding null-day with the same month value $m$.

Next, we propose two simple formulas which convert the century and year items into two integers. Recall that $c$ denotes the century item and $y\in[0,99]$ denotes the year item. We multiple the residue of $c$ modulo $4$ by $-2$ and denote the result by $w_1$. It is readily seen that
\begin{eqnarray}
  w_1(c)=\begin{cases}
    0,&\ c=4k;\\
    -2,&\ c=4k+1;\\
    -4,&\ c=4k+2;\\
    -6,&\ c=4k+3.
  \end{cases}
\end{eqnarray}
Our second formula is given as
\begin{equation}
  w_2(y)=\lfloor{5y\over4}\rfloor,
\end{equation}
where $\lfloor x\rfloor$ denotes the largest integer no more than $x$.
Let $y=10y_1+y_0$ with $y_0,y_1\in[0,9]$ being the ones and tens digits of $y$ respectively. The above formula can be modified as
\begin{equation}
  w_2(y)=\lfloor{25y_1\over2}+{5y_0\over4}\rfloor\equiv
y_0-y_1+\lfloor{y_0\over4}-{y_1\over2}\rfloor \mod 7.
\end{equation}

Finally, we compute $w(m,d,c,y)=w_0(m,d)+w_1(c)+w_2(y)$ modulo $7$ to obtain the day of the week for any given date $(m,d,c,y)$. If the given date is in January or February of a leap year, then we should subtract the sum by $1$, namely, $w(m,d,c,y)=w_0(m,d)+w_1(c)+w_2(y)-1$ modulo $7$.

\section{Examples}
We use two examples to illustrate our null-days algorithm.

Example 1. For the date March 26, 2014, we observe $m=3,~d=26,~c=20,~y=14$. Moreover, $y_0=4$ and $y_1=1$.
First, we obtain from the null-day $3/5$ that $w_0=21\equiv0\mod7$.
Next, we have $w_1=0$ since $c\equiv0\mod4$.
Furthermore, $w_2=4-1+\lfloor{4/4}-{1/2}\rfloor=3$.
Adding the above three numbers gives $w=w_0+w_1+w_2=3$, which implies that March 26, 2014 is a Wednesday.

Example 2. For the date February 10, 1984, we observe $m=2,~d=10,~c=19,~y=84$. Moreover, $y_0=4$ and $y_1=8$.
First, we obtain from the null-day $2/12$ that $w_0=-2\equiv5\mod7$.
Next, we have $w_1=-6$ since $c\equiv3\mod4$.
Furthermore, $w_2=4-8+\lfloor{4/4}-{8/2}\rfloor=-7\equiv0\mod7$.
Since $1984$ is a leap year and $m=2$, we add the above three numbers and subtract the sum by one to obtain $w=w_0+w_1+w_2-1=-2\equiv5\mod7$, which implies that February 10, 1984 is a Friday.



\begin{thebibliography}{1}

\bibitem{Ca87}
Lewis Carroll, To Find the Day of the Week for Any Given Date, {\it Nature}, {\bf 35} (1887), 517--517.

\bibitem{Co73}
John Horton Conway, Tomorrow is the Day After Doomsday, {\it Eureka}, {\bf 36} (1973), 28--31.

\end{thebibliography}
\end{document}